# ON MAXIMUM LIKELIHOOD ESTIMATION OF THE EXTREME VALUE INDEX

By Holger Drees[1], Ana Ferreira[2] and Laurens de Haan[2]

*University of Hamburg, EURANDOM and Technical University Lisbon, and Erasmus University Rotterdam*

We prove asymptotic normality of the so-called maximum likelihood estimator of the extreme value index.

**1. Introduction.** Let $X_1, X_2, \ldots$ be independent and identically distributed (i.i.d.) random variables (r.v.'s) from some unknown distribution function (d.f.) $F$. Denote the upper endpoint of $F$ by $x^*$, where $x^* = \sup\{x : F(x) < 1\} \leq \infty$, and let

$$(1) \qquad F_t(x) = P(X \leq t + x | X > t) = \frac{F(t+x) - F(t)}{1 - F(t)},$$

with $1 - F(t) > 0$, $t < x^*$ and $x > 0$, be the conditional d.f. of $X - t$ given $X > t$. Then it is well known [see Balkema and de Haan (1974) and Pickands (1975)] that up to scale and location transformations the generalized Pareto d.f. given by

$$(2) \qquad H_\gamma(x) = 1 - (1 + \gamma x)^{-1/\gamma},$$

$x > 0$ if $\gamma \geq 0$ and $0 < x < -1/\gamma$ if $\gamma < 0$ [for $\gamma = 0$ read $(1 + \gamma x)^{-1/\gamma}$ as $\exp(-x)$], can provide a good approximation to the conditional probabilities (1). More precisely, it has been proved that there exists a normalizing function $\sigma(t) > 0$, such that

$$\lim_{t \to x^*} F_t(x\sigma(t)) \to H_\gamma(x)$$

Received November 2002; revised June 2003.

[1]Supported by Netherlands Organization for Scientific Research through the Netherlands Mathematical Research Foundation and by the Heisenberg program of the DFG.

[2]Supported in part by POCTI/FCT/FEDER.

for all $x$, or equivalently

(3) $$\lim_{t \to x^*} \sup_{0 < x < x^* - t} |F_t(x) - H_\gamma(x/\sigma(t))| = 0,$$

if and only if $F$ is in the maximum domain of attraction of the corresponding extreme value d.f. $G_\gamma(x) = \exp(-(1+\gamma x)^{-1/\gamma})$ [Gnedenko (1943)], commonly denoted by $F \in D(G_\gamma)$. The parameter $\gamma \in \mathbb{R}$ is the extreme value index and is the same in both $H_\gamma$ and $G_\gamma$ approximations.

Under this set-up, it turns out that a major issue for estimating extreme events is the estimation of the extreme value index $\gamma$. A variety of procedures to estimate $\gamma$ are now available in the literature [e.g., Hill (1975), Dekkers, Einmahl and de Haan (1989) and Smith (1987)], although there are still open problems. Quite often the accuracy of these estimators relies heavily on the choice of some threshold, but it is not our aim here to address this type of optimality questions.

Instead, in this paper we present a relatively simple direct proof of the asymptotic normality of the maximum likelihood estimators (m.l.e.'s) of $\gamma$ and $\sigma$. It is based on some recent approximations to the tail empirical quantile function established by Drees (1998). Proofs of the asymptotic normality of the m.l.e.'s of $\gamma$ and $\sigma$ were given by Smith (1987), and also by Drees (1998) in the case $\gamma > 0$. Nonetheless we consider some proofs not easily understandable. Moreover, some of the conditions used in the aforementioned papers are unnecessarily restrictive.

For an i.i.d. sample of size $n$, let $X_{1,n} \leq X_{2,n} \leq \cdots \leq X_{n,n}$ be the ascending order statistics. In view of (3) we can expect that observations above some high threshold are approximately generalized Pareto. This motivates that inferences on $\gamma$ should be based on some set of high order statistics, say $(X_{n-k,n}, X_{n-k+1,n}, \ldots, X_{n,n})$, or, equivalently, on

$$Y_0 = X_{n-k,n},$$
$$Y_1 = X_{n-k+1,n} - X_{n-k,n},$$
$$\vdots$$
$$Y_k = X_{n,n} - X_{n-k,n},$$

where in the asymptotic setting $k = k_n$ is an intermediate sequence, that is, $k_n \to \infty$ and $k_n/n \to 0$ as $n \to \infty$. Since it is plausible that asymptotically the information contained in $Y_0$ is negligible as $k_n \to \infty$, we apply a conditional likelihood approach [see, e.g., Cox and Hinkley (1974), page 17] in that we consider the conditional distribution of the $(Y_1, \ldots, Y_k)$ given $Y_0 = y_0$. According to Theorem 2.4.1 of Arnold, Balakrishnan and Nagaraja (1992), it equals the distribution of the order statistics $(Y^*_{1,k}, \ldots, Y^*_{k,k})$ of an i.i.d. sample $(Y^*_1, \ldots, Y^*_k)$ with common distribution $F_{y_0}$ defined by (1). Hence, in view of convergence (3), the conditional distribution of the $(Y_1, \ldots, Y_k)$ given $Y_0 = y_0$ can be approximated by the distribution of an ordered sample of $k$ i.i.d. generalized Pareto random variables with d.f. $x \mapsto H_\gamma(x/\sigma)$. This suggests to estimate the unknown parameters $\gamma$ and $\sigma$ by a maximum likelihood estimator in the approximating generalized Pareto model; that is, given the sample $(x_1, \ldots, x_n)$ [or rather the largest observations $(x_{n-k,n}, \ldots, x_{n,n})$], we maximize $\prod_{i=1}^{k} h_{\gamma,\sigma}(y_i)$ with $y_i = x_{n-i+1,n} - x_{n-k,n}$, $1 \le i \le k$, and $h_{\gamma,\sigma}(y) = \partial H_\gamma(y/\sigma)/\partial y$.

Note that this approximative conditional likelihood function tends to $\infty$ if $\gamma < -1$ and $\gamma/\sigma \downarrow -1/(x_{n,n} - x_{n-k,n})$, and so a maximum over the full range of possible values for $(\gamma, \sigma)$ does not exist. Since, moreover, the maximum likelihood estimator behaves irregularly if $\gamma \le -1/2$, we look for a maximum of the approximative likelihood function only in the region $(\gamma, \sigma) \in (-1/2, \infty) \times (0, \infty)$.

The likelihood equations are then given in terms of the partial derivatives

$$\frac{\partial \log h_{\gamma,\sigma}(y)}{\partial \gamma} = \frac{1}{\gamma^2} \log\left(1 + \frac{\gamma}{\sigma} y\right) - \left(\frac{1}{\gamma} + 1\right) \frac{y/\sigma}{1 + (\gamma/\sigma) y},$$

$$\frac{\partial \log h_{\gamma,\sigma}(y)}{\partial \sigma} = -\frac{1}{\sigma} - \left(\frac{1}{\gamma} + 1\right) \frac{-(\gamma/\sigma^2) y}{1 + (\gamma/\sigma) y},$$

where for $\gamma = 0$ these terms should be interpreted as

$$\left.\frac{\partial \log h_{\gamma,\sigma}(y)}{\partial \gamma}\right|_{\gamma=0} = \frac{1}{2}\left(\frac{y}{\sigma}\right)^2 - \frac{y}{\sigma},$$

$$\left.\frac{\partial \log h_{\gamma,\sigma}(y)}{\partial \sigma}\right|_{\gamma=0} = -\frac{1}{\sigma} + \frac{y}{\sigma^2}.$$

The resulting likelihood equations in terms of the excesses $X_{n-i+1,n} - X_{n-k,n}$ are as follows:

(4)
$$\sum_{i=1}^{k} \frac{1}{\gamma^2} \log\left(1 + \frac{\gamma}{\sigma}(X_{n-i+1,n} - X_{n-k,n})\right)$$
$$- \left(\frac{1}{\gamma} + 1\right) \frac{(1/\sigma)(X_{n-i+1,n} - X_{n-k,n})}{1 + (\gamma/\sigma)(X_{n-i+1,n} - X_{n-k,n})} = 0,$$
$$\sum_{i=1}^{k} \left(\frac{1}{\gamma} + 1\right) \frac{(\gamma/\sigma)(X_{n-i+1,n} - X_{n-k,n})}{1 + (\gamma/\sigma)(X_{n-i+1,n} - X_{n-k,n})} = k$$

(with a similar interpretation when $\gamma = 0$), which for $\gamma \ne 0$ can be simplified to

$$\frac{1}{k} \sum_{i=1}^{k} \log\left(1 + \frac{\gamma}{\sigma}(X_{n-i+1,n} - X_{n-k,n})\right) = \gamma,$$

$$\frac{1}{k} \sum_{i=1}^{k} \frac{1}{1 + (\gamma/\sigma)(X_{n-i+1,n} - X_{n-k,n})} = \frac{1}{\gamma + 1},$$

with $(\gamma, \sigma) \in (-1/2, \infty) \times (0, \infty)$. The numerical problem to find a solution of these equations which maximizes the approximative likelihood was discussed by Grimshaw (1993).

From the above reasoning it follows that the m.l.e. of $\gamma$ is shift and scale invariant, and the m.l.e. of $\sigma$ is shift invariant and scale equivariant.

Next we sketch the proof of the asymptotic normality. Under standard second-order conditions [see (7)] we have for an intermediate sequence $k_n \in \mathbb{N}$

$$\text{(5)} \quad \left( \frac{Q_n(t) - Q_n(1)}{\tilde{a}(k_n/n)} \right)_{t \in [0,1]} = \left( \frac{t^{-\gamma_0} - 1}{\gamma_0} + k_n^{-1/2} Y_n(t) \right)_{t \in [0,1]},$$

where $(Q_n(t))_{t \in [0,1]}$ is a distributionally equivalent version of the process $(X_{n-[k_n t], n})_{t \in [0,1]}$, $(Y_n(t))_{t \in [0,1]}$ is an asymptotically Gaussian process of known mean and covariance function (Lemma 3.1), $\gamma_0$ is the true parameter and $\tilde{a}$ is a suitably chosen positive function [see (16)]. Hence for all $t \in [0,1]$ and all $\gamma$ and $\sigma$,

$$\text{(6)} \quad t^{\gamma_0} \left( 1 + \frac{\gamma}{\tilde{\sigma}} \frac{Q_n(t) - Q_n(1)}{\tilde{a}(k_n/n)} \right) = 1 + \left( \frac{\gamma}{\tilde{\sigma}} - \gamma_0 \right) \frac{1 - t^{\gamma_0}}{\gamma_0} + t^{\gamma_0} \frac{\gamma}{\tilde{\sigma}} k_n^{-1/2} Y_n(t),$$

where $\tilde{\sigma} = \sigma / \tilde{a}(k_n/n)$. Now if the sequence of solutions $(\gamma, \tilde{\sigma})$ satisfies

$$\gamma - \gamma_0 = O_p(k_n^{-1/2}) \quad \text{and} \quad \tilde{\sigma} - 1 = O_p(k_n^{-1/2}),$$

one can prove, using a construction similar to (5), that

$$\inf_{1/(2k_n) \leq t \leq 1} t^{\gamma_0} \left( 1 + \frac{\gamma}{\tilde{\sigma}} \frac{Q_n(t) - Q_n(1)}{\tilde{a}(k_n/n)} \right)$$

is stochastically bounded away from zero (Lemma 3.2). This implies

$$\log \left( t^{\gamma_0} \left( 1 + \frac{\gamma}{\tilde{\sigma}} \frac{Q_n(t) - Q_n(1)}{\tilde{a}(k_n/n)} \right) \right)$$

$$= \log \left( 1 + \left( \frac{\gamma}{\tilde{\sigma}} - \gamma_0 \right) \frac{1 - t^{\gamma_0}}{\gamma_0} + t^{\gamma_0} \frac{\gamma}{\tilde{\sigma}} k_n^{-1/2} Y_n(t) \right)$$

$$= \left( \frac{\gamma}{\tilde{\sigma}} - \gamma_0 \right) \frac{1 - t^{\gamma_0}}{\gamma_0} + t^{\gamma_0} \frac{\gamma}{\tilde{\sigma}} k_n^{-1/2} Y_n(t) + o_p(k_n^{-1/2})$$

and

$$\frac{1}{1 + (\gamma/\tilde{\sigma})(Q_n(t) - Q_n(1))/\tilde{a}(k_n/n)}$$

$$= t^{\gamma_0} \left( 1 - \left( \frac{\gamma}{\tilde{\sigma}} - \gamma_0 \right) \frac{1 - t^{\gamma_0}}{\gamma_0} - t^{\gamma_0} \frac{\gamma}{\tilde{\sigma}} k_n^{-1/2} Y_n(t) + o_p(k_n^{-1/2}) \right),$$

where the $o_p$-term is uniform for $1/(2k_n) \leq t \leq 1$ (proof of Proposition 3.1).

Hence, up to a $o_p(k_n^{-1/2})$-term, (4) are equivalent to linear equations which can be solved readily. The proof in Case $\gamma_0 = 0$ requires longer expansions but is similar.

The precise statement about the asymptotic normality is given in Theorem 2.1. In Theorem 2.2 an equivalent explicit estimator is constructed in the case $\gamma_0 = 0$.

Throughout the paper, $F^{\leftarrow}$ denotes the generalized inverse of $F$, $\xrightarrow{d}$ convergence in distribution and $\xrightarrow{p}$ convergence in probability.

**2. Asymptotic normality of the maximum likelihood estimators.** Assume that there exist measurable, locally bounded functions $a$, $\Phi:(0,1) \to (0,\infty)$ and $\Psi:(0,\infty) \to \mathbb{R}$ such that

$$(7) \qquad \lim_{t \downarrow 0} \frac{(F^{\leftarrow}(1-tx) - F^{\leftarrow}(1-t))/a(t) - (x^{-\gamma_0} - 1)/\gamma_0}{\Phi(t)} = \Psi(x),$$

for some $\gamma_0 > -1/2$, for all $t \in (0,1)$ and $x > 0$, where $x \mapsto \Psi(x)/(x^{-\gamma_0} - 1)$ is not constant, $\Phi(t)$ not changing sign eventually and $\Phi(t) \to 0$ as $t \downarrow 0$. Then, according to de Haan and Stadtmüller (1996), $|\Phi|$ is $-\rho$-varying at 0 for some $\rho \leq 0$, that is, $\lim_{t \downarrow 0} \Phi(tx)/\Phi(t) = x^{-\rho}$ for all $x > 0$, and

$$(8) \qquad \Psi(x) = \begin{cases} (x^{-(\gamma_0+\rho)} - 1)/(\gamma_0 + \rho), & \rho < 0, \\ -x^{-\gamma_0} \log(x)/\gamma_0, & \gamma_0 \neq \rho = 0, \\ \log^2(x), & \gamma_0 = \rho = 0, \end{cases}$$

provided that the normalizing function $a$ and the function $\Phi$ are chosen suitably. Condition (7) is a second-order refinement of $F \in D(G_{\gamma_0})$. Still, it is a quite general condition, satisfied for all usual distributions satisfying the max-domain of attraction condition.

We assume throughout that $k_n$ is an intermediate sequence, that is, $k_n \to \infty$ and $k_n/n \to 0$ as $n \to \infty$.

THEOREM 2.1. *Assume condition* (7) *for some* $\gamma_0 > -1/2$ *and that the intermediate sequence* $k_n$ *satisfies*

$$(9) \qquad \Phi(k_n/n) = O(k_n^{-1/2}).$$

*Then the system of likelihood equations* (4) *has a sequence of solutions* $(\hat{\gamma}_n, \hat{\sigma}_n)$ *that verifies*

$$(10) \quad \begin{aligned} & k_n^{1/2}(\hat{\gamma}_n - \gamma_0) \\ & \quad - \frac{(\gamma_0 + 1)^2}{\gamma_0} k_n^{1/2} \Phi\left(\frac{k_n}{n}\right) \int_0^1 (t^{\gamma_0} - (2\gamma_0 + 1)t^{2\gamma_0}) \Psi(t) \, dt \\ & \quad \xrightarrow{d} \frac{(\gamma_0 + 1)^2}{\gamma_0} \int_0^1 (t^{\gamma_0} - (2\gamma_0 + 1)t^{2\gamma_0})(W(1) - t^{-(\gamma_0+1)}W(t)) \, dt, \end{aligned}$$

$$k_n^{1/2}\left(\frac{\hat{\sigma}_n}{a(k_n/n)} - 1\right)$$
$$- \frac{\gamma_0 + 1}{\gamma_0} k_n^{1/2} \Phi\left(\frac{k_n}{n}\right) \int_0^1 ((\gamma_0+1)(2\gamma_0+1)t^{2\gamma_0} - t^{\gamma_0})\Psi(t)\,dt$$
$$\stackrel{d}{\to} \frac{\gamma_0+1}{\gamma_0} \int_0^1 ((\gamma_0+1)(2\gamma_0+1)t^{2\gamma_0} - t^{\gamma_0})(W(1) - t^{-(\gamma_0+1)}W(t))\,dt,$$

(11)

*as* $n \to \infty$, *and the convergence holds jointly with the same standard Brownian motion* $W$. *For* $\gamma_0 = 0$ *these equations should be interpreted as their limits when* $\gamma_0 \to 0$; *that is*,

(12)
$$k_n^{1/2}\hat{\gamma}_n + k_n^{1/2}\Phi\left(\frac{k_n}{n}\right)\int_0^1 (2+\log t)\Psi(t)\,dt$$
$$\stackrel{d}{\to} -\int_0^1 (2+\log t)(W(1) - t^{-1}W(t))\,dt,$$

(13)
$$k_n^{1/2}\left(\frac{\hat{\sigma}_n}{a(k_n/n)} - 1\right) - k_n^{1/2}\Phi\left(\frac{k_n}{n}\right)\int_0^1 (3+\log t)\Psi(t)\,dt$$
$$\stackrel{d}{\to} \int_0^1 (3+\log t)(W(1) - t^{-1}W(t))\,dt.$$

*Moreover, any sequence of solutions* $(\hat{\gamma}_n^*, \hat{\sigma}_n^*)$ *which is not of the type* (10)–(13) *must satisfy* $k_n^{1/2}|\hat{\gamma}_n^* - \gamma_0| \stackrel{p}{\to} \infty$ *or* $k_n^{1/2}|\hat{\sigma}_n^*/a(k_n/n) - 1| \stackrel{p}{\to} \infty$.

REMARK 2.1. Condition (9) is satisfied if $k_n \to \infty$ not too fast. The bias term $(\gamma_0+1)^2 k_n^{1/2}\Phi(k_n/n)\int_0^1 (t^{\gamma_0} - (2\gamma_0+1)t^{2\gamma_0})\Psi(t)\,dt/\gamma_0$ in (10) vanishes if $k_n^{1/2}\Phi(k_n/n) \to 0$. A similar remark applies to (11)–(13).

REMARK 2.2. Note that the likelihood equations are satisfied with $\gamma = 0$ if and only if

$$\frac{1}{2k}\sum_{i=1}^k (X_{n-i+1,n} - X_{n-k,n})^2 = \left(\frac{1}{k}\sum_{i=1}^k (X_{n-i+1,n} - X_{n-k,n})\right)^2$$

and $\sigma = \sum_{i=1}^k (X_{n-i+1,n} - X_{n-k,n})/k$. Hence, the m.l.e. for $\gamma$ will a.s. not be equal to 0 if, for example, $F$ possesses a density.

COROLLARY 2.1. *Under the conditions of Theorem* 2.1 *and if*

(14)
$$k_n^{1/2}\Phi\left(\frac{k_n}{n}\right) \to \lambda \in \mathbb{R},$$

*the solutions* (10)–(13) *verify*

(15)
$$k_n^{1/2}\left[\begin{array}{c}\hat{\gamma}_n - \gamma_0 \\ \hat{\sigma}_n/a(k_n/n) - 1\end{array}\right] \stackrel{d}{\to} N(\lambda\mu, \Sigma),$$

*where $N$ denotes the bivariate normal distribution, $\mu$ equals*

$$\begin{cases} \left[\dfrac{\rho(\gamma_0+1)}{(1-\rho)(\gamma_0-\rho+1)}, \dfrac{1-2\rho+\gamma_0-\rho\gamma_0}{(1-\rho)(\gamma_0-\rho+1)}\right]^T, & \text{if } \rho < 0, \\ [1, \gamma_0^{-1}]^T, & \text{if } \gamma_0 \neq \rho = 0, \\ [2, 0]^T, & \text{if } \gamma_0 = \rho = 0, \end{cases}$$

*and*

$$\Sigma = \begin{bmatrix} (1+\gamma_0)^2 & -(1+\gamma_0) \\ -(1+\gamma_0) & 2+2\gamma_0+\gamma_0^2 \end{bmatrix}.$$

REMARK 2.3. Smith (1987) examined a slightly different version of the m.l.e. that is based on the excesses over a *deterministic* threshold $u = u_n$ instead of the excesses over the random threshold $X_{n-k_n,n}$. For the comparison of Smith's results with Corollary 2.1, we focus on the case $\gamma_0 \neq 0$, $\rho < 0$ and $\lambda = 0$, when there is no asymptotic bias, since in the other cases the more restrictive second-order conditions used by Smith are not directly comparable to our setting.

Let $K$ denote the (random) number of exceedances over the threshold $u$ and let

$$\sigma_n = \begin{cases} \gamma_0 u, & \gamma_0 > 0, \\ |\gamma_0|(F^{\leftarrow}(1) - u), & \gamma_0 < 0. \end{cases}$$

Then it was shown that the standardized m.l.e.'s $K^{1/2}(\hat{\gamma}_n - \gamma_0, \hat{\sigma}_n/\sigma_n - 1)$ based on the exceedances $X_i - u$ converge to a centered bivariate normal distribution with covariance matrix

$$\begin{pmatrix} (1+\gamma_0)^2 & -(1+\gamma_0) \\ -(1+\gamma_0) & 2(1+\gamma_0) \end{pmatrix}.$$

At first glance, it seems peculiar that we obtain a different asymptotic variance for the scale estimator in Corollary 2.1, namely $2(1+\gamma_0)+\gamma_0^2$. However, the following heuristic reasoning shows that in fact the increase in the variance is due to the slightly different standardization.

To make the results about the asymptotic behavior comparable, in our setting one has to condition at the event $X_{n-k_n,n} = u$. Then Smith's result claims that conditionally $\hat{\sigma}_n = \sigma_n(1 + k_n^{-1/2} Z_n)$ for some asymptotically centered normal r.v. $Z_n$ with asymptotic variance $2(1+\gamma_0)$. Hence conditionally at $X_{n-k_n,n} = u$,

$$k_n^{1/2}\left(\frac{\hat{\sigma}_n}{a(k_n/n)} - 1\right) = \frac{\sigma_n}{a(k_n/n)} Z_n + k_n^{1/2}\left(\frac{\sigma_n}{a(k_n/n)} - 1\right).$$

Because, in the restrictive setting considered here, $a(k_n/n) = \gamma_0 F^{\leftarrow}(1 - k_n/n)$ for $\gamma_0 > 0$ and $a(k_n/n) = |\gamma_0|(F^{\leftarrow}(1) - F^{\leftarrow}(1 - k_n/n))$ for $\gamma_0 < 0$, unconditionally (i.e., when $u$ is replaced with $X_{n-k_n,n}$) $\sigma_n/a(k_n/n) \to 1$ in probability, so that the first term tends to a normal random variable with variance $2(\gamma_0 + 1)$. According to the approximation of the tail empirical quantile function [cf. (18)], unconditionally the second term converges to $\gamma_0 W(1)$. Since asymptotically $X_{n-k_n,n}$ and the excesses $X_{n-i+1,n} - X_{n-k_n,n}, 1 \leq i \leq k_n$, are independent, so are $Z_n$ and $W(1)$. Hence the two variances $2(\gamma_0 + 1)$ and $\gamma_0^2$ add up, leading to the variance given in Corollary 2.1.

We now show that if $\gamma_0 = 0$, the m.l.e.'s are asymptotically equivalent in some sense to explicit estimators. Define

$$m_n^{(j)} = \frac{1}{k_n} \sum_{i=1}^{k_n} (X_{n-i+1,n} - X_{n-k_n,n})^j, \qquad j = 1, 2,$$

$$\hat{\gamma}_* = 1 - \frac{1}{2}\left(1 - \frac{(m_n^{(1)})^2}{m_n^{(2)}}\right)^{-1}$$

and

$$\hat{a}_*\left(\frac{k_n}{n}\right) = \frac{2(m_n^{(1)})^3}{m_n^{(2)}}.$$

It can be shown, using Corollary 3.1, that these estimators are consistent and asymptotically normal if $\gamma < 1/2$. Let $(\hat{\gamma}_{\text{MLE}}, \hat{\sigma}_{\text{MLE}})$ be a sequence of solutions of (4) as described in Theorem 2.1.

THEOREM 2.2. *If $F$ is in the class of distributions that satisfy (7) with $\gamma_0 = 0$ and if (9) holds, then*

$$k_n^{1/2}(\hat{\gamma}_* - \hat{\gamma}_{\text{MLE}}) \xrightarrow{p} 0$$

*and*

$$k_n^{1/2}\left(\frac{\hat{a}_*(k_n/n) - \hat{\sigma}_{\text{MLE}}}{a(k_n/n)}\right) \xrightarrow{p} 0.$$

REMARK 2.4. If, in addition, (7) holds with $\rho < 0$, $\sup\{x|F(x) < 1\} > 0$ and $k_n = o(\log^2 n)$, then we have an analogous result for the moment estimator introduced by Dekkers, Einmahl and de Haan (1989):

$$k_n^{1/2}(\hat{\gamma}_{\text{MOM}} - \hat{\gamma}_{\text{MLE}}) \xrightarrow{p} 0,$$

where

$$\hat{\gamma}_{\text{MOM}} = M_n^{(1)} + 1 - \frac{1}{2}\left(1 - \frac{(M_n^{(1)})^2}{M_n^{(2)}}\right)^{-1},$$

with $M_n^{(j)} = \frac{1}{k_n} \sum_{i=0}^{k_n-1} (\log X_{n-i,n} - \log X_{n-k_n,n})^j$, $j = 1, 2$. A similar statement holds for the scale estimator

$$\hat{a}_{**}\left(\frac{k_n}{n}\right) = \frac{2(M_n^{(1)})^3}{M_n^{(2)}}.$$

The condition $k_n = o(\log^2 n)$ ensures that the bias vanishes asymptotically. We prove this remark in Section 3.

**3. Proofs.** Given (7) with $\gamma_0 > -1/2$ and (9), from Theorem 2.1 in Drees (1998) one can find a probability space and define on that space a Brownian motion $W$ and a sequence of stochastic processes $Q_n$ such that (i) for each $n$, $(Q_n(t))_{t \in [0,1]} \stackrel{d}{=} (X_{n-[k_n t],n})_{t \in [0,1]}$, and (ii) there exist functions $\tilde{a}(k_n/n) = a(k_n/n)(1 + o(\Phi(k_n/n)))$ and $\tilde{\Phi}(k_n/n) \sim \Phi(k_n/n)$ such that, for all $\varepsilon > 0$,

(16)
$$\sup_{t \in [0,1]} t^{\gamma_0 + 1/2 + \varepsilon} \left| \frac{Q_n(t) - F^{\leftarrow}(1 - k_n/n)}{\tilde{a}(k_n/n)} - \left( \frac{t^{-\gamma_0} - 1}{\gamma_0} - t^{-(\gamma_0+1)} \frac{W(k_n t)}{k_n} + \tilde{\Phi}\left(\frac{k_n}{n}\right) \Psi(t) \right) \right|$$
$$= o_p(k_n^{-1/2}) + o_p\left( \tilde{\Phi}\left(\frac{k_n}{n}\right) \right) \quad \text{as } n \to \infty.$$

A similar expansion is also valid for $\gamma_0 \leq -1/2$ when $F^{\leftarrow}(1 - k_n/n)$ is replaced with a suitable random variable.

Define

(17) $$Y_n(t) = k_n^{1/2} \left( \frac{Q_n(t) - Q_n(1)}{\tilde{a}(k_n/n)} - \frac{t^{-\gamma_0} - 1}{\gamma_0} \right)$$

[read $(t^{-\gamma_0} - 1)/\gamma_0$ as $-\log t$, when $\gamma_0 = 0$]. Hence we have the following lemma.

LEMMA 3.1. *Suppose* (7) *holds and that the intermediate sequence* $k_n$ *satisfies* (9). *Then, for all* $\varepsilon > 0$,

(18) $$Y_n(t) = W_n(1) - t^{-(\gamma_0+1)} W_n(t) + k_n^{1/2} \tilde{\Phi}\left(\frac{k_n}{n}\right) \Psi(t) + o_p(1) t^{-(\gamma_0+1/2+\varepsilon)}$$

*as* $n \to \infty$, *where* $W_n(t) = k_n^{-1/2} W(k_n t)$ *is a standard Brownian motion and the* $o_p$-*term is uniform for* $t \in [0, 1]$.

From this lemma the following corollary follows easily.

COROLLARY 3.1. *Under the conditions of Lemma* 3.1, *for all* $\varepsilon > 0$,

(19) $$Y_n(t) = O_p(1) t^{-(\gamma_0+1/2+\varepsilon)}$$

*as* $n \to \infty$, *where the* $O_p$-*term is uniform for* $t \in [0, 1]$.

Given the previous results, to prove Theorem 2.1 it is sufficient to consider the likelihood equations with $(X_{n-[k_n t],n} - X_{n-k_n,n})$ replaced by $Q_n(t) - Q_n(1)$, $t \in [0,1]$. It is convenient to reparametrize the equations in terms of $(\gamma, \tilde{\sigma}) = (\gamma, \sigma/\tilde{a}(k_n/n))$. Then we have the equations

$$
\begin{aligned}
&\int_0^1 \left( \frac{1}{\gamma^2} \log\left(1 + \frac{\gamma}{\tilde{\sigma}} \frac{Q_n(t) - Q_n(1)}{\tilde{a}(k_n/n)}\right) \right. \\
&\quad \left. - \left(\frac{1}{\gamma} + 1\right) \frac{(1/\tilde{\sigma})(Q_n(t) - Q_n(1))/\tilde{a}(k_n/n)}{1 + (\gamma/\tilde{\sigma})(Q_n(t) - Q_n(1))/\tilde{a}(k_n/n)} \right) dt = 0, \\
&\int_0^1 \left(\frac{1}{\gamma} + 1\right) \frac{(\gamma/\tilde{\sigma})(Q_n(t) - Q_n(1))/\tilde{a}(k_n/n)}{1 + (\gamma/\tilde{\sigma})(Q_n(t) - Q_n(1))/\tilde{a}(k_n/n)} dt = 1.
\end{aligned}
\tag{20}
$$

LEMMA 3.2. *Assume conditions* (7) *and* (9). *Let* $(\gamma, \tilde{\sigma}) = (\gamma_n, \tilde{\sigma}_n)$ *be such that*

$$|\gamma/\tilde{\sigma} - \gamma_0| = O_p(k_n^{-1/2}). \tag{21}$$

*Then, if* $-1/2 < \gamma_0 < 0$ *or* $\gamma_0 > 0$,

$$P\left(1 + \frac{\gamma}{\tilde{\sigma}} \frac{Q_n(t) - Q_n(1)}{\tilde{a}(k_n/n)} \geq C_n t^{-\gamma_0}, t \in \left[\frac{1}{2k_n}, 1\right]\right) \to 1, \qquad n \to \infty, \tag{22}$$

*for some r.v.'s* $C_n > 0$ *such that* $1/C_n = O_P(1)$. *If* $\gamma_0 = 0$,

$$P\left(1 + \frac{\gamma}{\tilde{\sigma}} \frac{Q_n(t) - Q_n(1)}{\tilde{a}(k_n/n)} \geq \frac{1}{2}, t \in \left[\frac{1}{2k_n}, 1\right]\right) \to 1, \qquad n \to \infty, \tag{23}$$

*and*

$$\sup_{t \in [0,1]} \frac{Q_n(t) - Q_n(1)}{\tilde{a}(k_n/n)} = O_p(\log k_n), \qquad n \to \infty. \tag{24}$$

PROOF. It suffices to prove the assertions with $Q_n(t)$ replaced with $X_{n-[k_n t],n}$. Without loss of generality, one may assume $X_{i,n} = F^{\leftarrow}(1 - U_{i,n})$ for uniform order statistics $U_{i,n}$ since $(X_{n-[k_n t],n})_{t \in [0,1]} \stackrel{d}{=} (F^{\leftarrow}(1 - U_{[k_n t]+1,n}))_{t \in [0,1]} \stackrel{d}{=} (Q_n(t))_{t \in [0,1]}$.

Note that, by Shorack and Wellner [(1986), Chapter 10, Section 3, page 416, inequality 2],

$$\sup_{1/(2k_n) \leq t \leq 1} \frac{nU_{[k_n t]+1,n}}{k_n t} = O_P(1), \qquad \sup_{0 \leq t \leq 1} \frac{k_n t}{nU_{[k_n t]+1,n}} = O_P(1), \tag{25}$$

as $n \to \infty$. Also note that (7) implies, for some functions $\tilde{a}(s) \sim a(s)$ and $\tilde{\Phi}(s) \sim \Phi(s)$, $s \downarrow 0$, for all $x_0 > 0$ and $\varepsilon > 0$,

$$\lim_{s \downarrow 0} \sup_{0 < x \leq x_0} x^{\gamma_0 + \varepsilon} \left| \frac{(F^{\leftarrow}(1 - sx) - F^{\leftarrow}(1 - s))/\tilde{a}(s) - (x^{-\gamma_0} - 1)/\gamma_0}{\tilde{\Phi}(s)} \right.$$

$$\left. - \Psi(x) \right| = 0$$

[Drees (1998), Lemma 2.1]. Combining these two results, we obtain

$$\sup_{t \in [1/(2k_n), 1]} t^{\gamma_0 + \varepsilon} \left| \left( \frac{F^\leftarrow(1 - U_{[k_n t]+1, n}) - F^\leftarrow(1 - k_n/n)}{\tilde{a}(k_n/n)} \right. \right.$$
$$\left. \left. - \frac{((n/k_n) U_{[k_n t]+1, n})^{-\gamma_0} - 1}{\gamma_0} \right) \right.$$
$$\left. \times \tilde{\Phi}\left(\frac{k_n}{n}\right)^{-1} - \Psi\left(\frac{n}{k_n} U_{[k_n t]+1, n}\right) \right| = o_p(1). \tag{26}$$

Next use this approximation simultaneously for $t \in [1/(2k_n), 1]$ and $t = 1$. In view of the special construction $X_{n-[k_n t], n} = F^\leftarrow(1 - U_{[k_n t]+1, n})$, we then have, for $-1/2 < \gamma_0 < 0$ or $\gamma_0 > 0$,

$$\frac{X_{n-[k_n t], n} - X_{n-k_n, n}}{\tilde{a}(k_n/n)}$$
$$= \frac{F^\leftarrow(1 - U_{[k_n t]+1, n}) - F^\leftarrow(1 - U_{k_n+1, n})}{\tilde{a}(k_n/n)}$$
$$= \frac{1}{\gamma_0}\left(\frac{n}{k_n} U_{[k_n t]+1, n}\right)^{-\gamma_0} - \frac{1}{\gamma_0}\left(\frac{n}{k_n} U_{k_n+1, n}\right)^{-\gamma_0}$$
$$+ \tilde{\Phi}\left(\frac{k_n}{n}\right) \Psi\left(\frac{n}{k_n} U_{[k_n t]+1, n}\right) - \tilde{\Phi}\left(\frac{k_n}{n}\right) \Psi\left(\frac{n}{k_n} U_{k_n+1, n}\right)$$
$$+ o_p\left(t^{-(\gamma_0 + \varepsilon)} \tilde{\Phi}\left(\frac{k_n}{n}\right)\right).$$

Hence

$$1 + \frac{\gamma}{\tilde{\sigma}} \frac{X_{n-[k_n t], n} - X_{n-k_n, n}}{\tilde{a}(k_n/n)}$$
$$= \left(1 - \left(\frac{n}{k_n} U_{k_n+1, n}\right)^{-\gamma_0}\right) - \left(\frac{\gamma}{\tilde{\sigma}} - \gamma_0\right) \frac{1}{\gamma_0} \left(\frac{n}{k_n} U_{k_n+1, n}\right)^{-\gamma_0}$$
$$+ \frac{\gamma}{\tilde{\sigma}} \frac{1}{\gamma_0}\left(\frac{n}{k_n} U_{[k_n t]+1, n}\right)^{-\gamma_0} + \frac{\gamma}{\tilde{\sigma}} \tilde{\Phi}\left(\frac{k_n}{n}\right) \Psi\left(\frac{n}{k_n} U_{[k_n t]+1, n}\right)$$
$$- \frac{\gamma}{\tilde{\sigma}} \tilde{\Phi}\left(\frac{k_n}{n}\right) \Psi\left(\frac{n}{k_n} U_{k_n+1, n}\right) + o_p(t^{-(\gamma_0 + \varepsilon)} k_n^{-1/2})$$
$$= I + II + III + IV + V + VI.$$

By (25), $t^{\gamma_0} III$ is bounded away from zero uniformly for $t \in [(2k_n)^{-1}, 1]$. We will show that all the other terms tend to 0 uniformly when multiplied with $t_0^\gamma$, so that assertion (22) follows with $C_n := \inf_{t \in [(2k_n)^{-1}, 1]} t^{\gamma_0} III - \varepsilon_n$ for a suitable sequence $\varepsilon_n \downarrow 0$.

By the asymptotic normality of intermediate order statistics, part $I$ is $O_p(k_n^{-1/2})$. Hence $t^{\gamma_0} I = o_p(1)$, which is trivial if $\gamma_0 > 0$; for $-1/2 < \gamma_0 < 0$ note that $t^{\gamma_0} k_n^{-1/2} \leq 2^{-\gamma_0} k_n^{-\gamma_0-1/2} \to 0$ as $k_n \to \infty$. By (25) and assumption (21), part $II$ is $O_p(k_n^{-1/2})$ so that by the same arguments as above, $t^{\gamma_0} II = o_P(1)$.

Next note that $t^{\gamma_0} \Psi(t) = o(t^{-1/2})$ as $t \downarrow 0$. This combined with (9) and (25) gives that $t^{\gamma_0} IV$ and $t^{\gamma_0} V$ are $o_p(1)$. Finally, $t^{\gamma_0} VI = o_p(1)$, provided one chooses $\varepsilon < 1/2$.

Now consider the case $\gamma_0 = 0$. Since (26) is still valid when $\gamma_0 = 0$, with the obvious changes, we get

$$
\begin{aligned}
(27) \quad & \frac{X_{n-[k_n t],n} - X_{n-k_n,n}}{\tilde{a}(k_n/n)} \\
& = -\log\left(\frac{n}{k_n} U_{[k_n t]+1,n}\right) + \log\left(\frac{n}{k_n} U_{k_n+1,n}\right) \\
& \quad + \tilde{\Phi}\left(\frac{k_n}{n}\right) \Psi\left(\frac{n}{k_n} U_{[k_n t]+1,n}\right) - \tilde{\Phi}\left(\frac{k_n}{n}\right) \Psi\left(\frac{n}{k_n} U_{k_n+1,n}\right) \\
& \quad + o_p\left(t^{-\varepsilon} \tilde{\Phi}\left(\frac{k_n}{n}\right)\right).
\end{aligned}
$$

Hence

$$
\begin{aligned}
1 & + \frac{\gamma}{\tilde{\sigma}} \frac{X_{n-[k_n t],n} - X_{n-k_n,n}}{\tilde{a}(k_n/n)} \\
& = 1 - \frac{\gamma}{\tilde{\sigma}} \log t - \frac{\gamma}{\tilde{\sigma}} \log\left(\frac{n}{k_n t} U_{[k_n t]+1,n}\right) + \frac{\gamma}{\tilde{\sigma}} \log\left(\frac{n}{k_n} U_{k_n+1,n}\right) \\
& \quad + \frac{\gamma}{\tilde{\sigma}} \tilde{\Phi}\left(\frac{k_n}{n}\right) \Psi\left(\frac{n}{k_n} U_{[k_n t]+1,n}\right) - \frac{\gamma}{\tilde{\sigma}} \tilde{\Phi}\left(\frac{k_n}{n}\right) \Psi\left(\frac{n}{k_n} U_{k_n+1,n}\right) \\
& \quad + o_p\left(\frac{\gamma}{\tilde{\sigma}} t^{-\varepsilon} k_n^{-1/2}\right).
\end{aligned}
$$

Hence by (25) and assumptions (9) and (21), all the terms but the 1 in the last equality tend to 0 in probability uniformly for $t \in [1/(2k_n), 1]$ so that (23) is obvious.

Finally, to verify (24) just note that, for $t = 1/(2k_n)$, the expression (27) is of the order $O_p(\log k_n)$, provided $0 < \varepsilon < 1/2$. Since $X_{n-[k_n t],n} \leq X_{n,n}$ for all $t \in [0,1]$ the assertion (24) follows. $\square$

PROPOSITION 3.1.  *Assume conditions* (7) *and* (9). *Any solution* $(\gamma, \tilde{\sigma})$ *of* (20) *satisfying* (21) *and* $\log \tilde{\sigma} = O_P(1)$ *admits the approximation*

$$
\begin{aligned}
(28) \quad & k_n^{1/2}(\gamma - \gamma_0) - \frac{(\gamma_0+1)^2}{\gamma_0} \int_0^1 (t^{\gamma_0} - (2\gamma_0+1)t^{2\gamma_0}) Y_n(t)\, dt = o_p(1), \\
& k_n^{1/2}(\tilde{\sigma} - 1) - \frac{\gamma_0+1}{\gamma_0} \int_0^1 ((\gamma_0+1)(2\gamma_0+1)t^{2\gamma_0} - t^{\gamma_0}) Y_n(t)\, dt = o_p(1),
\end{aligned}
$$

as $n \to \infty$. *For $\gamma_0 = 0$ these equations should be interpreted as their limits for $\gamma_0 \to 0$, that is,*

(29)
$$k_n^{1/2}\gamma + \int_0^1 (2 + \log t)Y_n(t)\,dt = o_p(1),$$
$$k_n^{1/2}(\tilde{\sigma} - 1) - \int_0^1 (3 + \log t)Y_n(t)\,dt = o_p(1).$$

*Conversely, there exists a solution of* (20) *which satisfies* (28), *respectively* (29), *and hence also* (21).

REMARK 3.1. For $\gamma_0 \neq 0$ the condition on $\log \tilde{\sigma}$ is not needed.

PROOF OF PROPOSITION 3.1. We consider the cases $\gamma_0 > 0$, $-1/2 < \gamma_0 < 0$ and $\gamma_0 = 0$ separately.

CASE $\gamma_0 > 0$. In view of assumption (21), we may assume $\gamma \neq 0$. Hence, system (20) can be simplified to

(30)
$$\int_0^1 \log\left(1 + \frac{\gamma}{\tilde{\sigma}} \frac{Q_n(t) - Q_n(1)}{\tilde{a}(k_n/n)}\right) dt = \gamma,$$
$$\int_0^1 \frac{1}{1 + (\gamma/\tilde{\sigma})(Q_n(t) - Q_n(1))/\tilde{a}(k_n/n)}\,dt = \frac{1}{\gamma + 1}.$$

Next we will find expansions for the left-hand side of both equations.

Rewrite the first one as

$$\int_0^{(2k_n)^{-1}} \log\left(1 + \frac{\gamma}{\tilde{\sigma}} \frac{Q_n(t) - Q_n(1)}{\tilde{a}(k_n/n)}\right) dt + \int_{(2k_n)^{-1}}^1 \log(t^{-\gamma_0})\,dt$$
$$+ \int_{(2k_n)^{-1}}^1 \log\left(t^{\gamma_0}\left(1 + \frac{\gamma}{\tilde{\sigma}} \frac{Q_n(t) - Q_n(1)}{\tilde{a}(k_n/n)}\right)\right) dt$$
$$= I_1 + \gamma_0(1 - O(k_n^{-1}\log k_n)) + I_2.$$

First we prove that $I_1$ is negligible. Since $t \mapsto Q_n(t)$ is constant when $t \in [0, (2k_n)^{-1}]$, Lemma 3.2 implies that, with probability tending to 1,

(31) $1 + \dfrac{\gamma}{\tilde{\sigma}} \dfrac{Q_n(t) - Q_n(1)}{\tilde{a}(k_n/n)} = 1 + \dfrac{\gamma}{\tilde{\sigma}} \dfrac{Q_n(1/(2k_n)) - Q_n(1)}{\tilde{a}(k_n/n)} \geq (2k_n)^{\gamma_0} C_n$

for all $t \in [0, (2k_n)^{-1}]$ with $C_n$ stochastically bounded away from 0, so that $-I_1 \leq (2k_n)^{-1}O_P(\log k_n)$. On the other hand, from (17), (19) and (21),

$$1 + \frac{\gamma}{\tilde{\sigma}} \frac{Q_n(1/(2k_n)) - Q_n(1)}{\tilde{a}(k_n/n)}$$
$$= 1 + (\gamma_0 + O_P(k_n^{-1/2}))\left(\frac{(2k_n)^{\gamma_0} - 1}{\gamma_0} + O_P(k_n^{\gamma_0+\varepsilon})\right) = O_P(k_n^{\gamma_0+\varepsilon}).$$

Hence, it follows that $I_1 = o_p(k_n^{-1/2})$.

Next we turn to the main term $I_2$. We will apply the inequality $0 \leq x - \log(1+x) \leq x^2/(2(1 \wedge (1+x)))$, valid for all $x > -1$, to

$$\begin{aligned}
x &= t^{\gamma_0}\left(1 + \frac{\gamma}{\tilde{\sigma}}\frac{Q_n(t) - Q_n(1)}{\tilde{a}(k_n/n)}\right) - 1 \\
&= \left(\frac{\gamma}{\tilde{\sigma}} - \gamma_0\right)\frac{1 - t^{\gamma_0}}{\gamma_0} + \frac{\gamma}{\tilde{\sigma}}k_n^{-1/2}t^{\gamma_0}Y_n(t).
\end{aligned} \tag{32}$$

Then, from Lemma 3.2 it follows that $0 < 1/(1 \wedge (1+x)) \leq 1 \vee 1/C_n = O_P(1)$ with probability tending to 1. Moreover, note that relation (19) implies

$$\int_0^{(2k_n)^{-1}} t^{\gamma_0} Y_n(t)\, dt = O_P\left(\int_0^{(2k_n)^{-1}} t^{-1/2-\varepsilon}\, dt\right) = O_p((2k_n)^{-1/2+\varepsilon}) = o_p(1),$$

for $\varepsilon \in (0, 1/2)$. Hence from (17) and (19), as $n \to \infty$,

$$\begin{aligned}
I_2 &= \int_{(2k_n)^{-1}}^1 \left(t^{\gamma_0}\left(1 + \frac{\gamma}{\tilde{\sigma}}\frac{Q_n(t) - Q_n(1)}{\tilde{a}(k_n/n)}\right) - 1\right) dt \\
&\quad + O_P\left(\int_{(2k_n)^{-1}}^1 \left(t^{\gamma_0}\left(1 + \frac{\gamma}{\tilde{\sigma}}\frac{Q_n(t) - Q_n(1)}{\tilde{a}(k_n/n)}\right) - 1\right)^2 dt\right) \\
&= \int_{(2k_n)^{-1}}^1 \left(\left(\frac{\gamma}{\tilde{\sigma}} - \gamma_0\right)\frac{1 - t^{\gamma_0}}{\gamma_0} + \frac{\gamma}{\tilde{\sigma}}k_n^{-1/2}t^{\gamma_0}Y_n(t)\right) dt \\
&\quad + O_p\left(\int_{(2k_n)^{-1}}^1 \left(\left(\frac{\gamma}{\tilde{\sigma}} - \gamma_0\right)\frac{1 - t^{\gamma_0}}{\gamma_0} + \frac{\gamma}{\tilde{\sigma}}k_n^{-1/2}t^{\gamma_0}Y_n(t)\right)^2 dt\right) \\
&= \left(\left(\frac{\gamma}{\tilde{\sigma}} - \gamma_0\right)\frac{1}{\gamma_0 + 1} + O_p(k_n^{-1/2}(2k_n)^{-1})\right) \\
&\quad + \left(\frac{\gamma}{\tilde{\sigma}}k_n^{-1/2}\int_0^1 t^{\gamma_0}Y_n(t)\, dt + o_p(k_n^{-1/2})\right) \\
&\quad + O_p(k_n^{-1} + k_n^{-1}(2k_n)^{2\varepsilon} + k_n^{-1}(2k_n)^{-1/2+\varepsilon}) \\
&= \left(\frac{\gamma}{\tilde{\sigma}} - \gamma_0\right)\frac{1}{\gamma_0 + 1} + \frac{\gamma}{\tilde{\sigma}}k_n^{-1/2}\int_0^1 t^{\gamma_0}Y_n(t)\, dt + o_p(k_n^{-1/2}),
\end{aligned}$$

where for the last equality we took $\varepsilon < 1/4$. To sum up, we have proved that

$$\int_0^1 \log\left(1 + \frac{\gamma}{\tilde{\sigma}}\frac{Q_n(t) - Q_n(1)}{\tilde{a}(k_n/n)}\right) dt$$

$$= \gamma_0 + \left(\frac{\gamma}{\tilde{\sigma}} - \gamma_0\right)\frac{1}{\gamma_0 + 1} + \frac{\gamma}{\tilde{\sigma}}k_n^{-1/2}\int_0^1 t^{\gamma_0}Y_n(t)\, dt + o_p(k_n^{-1/2}).$$

This means that the first equation of (30) is equivalent to

$$\gamma = \gamma_0 + \left(\frac{\gamma}{\tilde{\sigma}} - \gamma_0\right)\frac{1}{\gamma_0 + 1}$$

$$+ \frac{\gamma}{\tilde{\sigma}} k_n^{-1/2} \int_0^1 t^{\gamma_0} Y_n(t)\,dt + o_p(k_n^{-1/2}).$$

Now we deal with the left-hand side of the second equation in (30). Applying the equality

$$\frac{1}{1+x} = 1 - x + \frac{x^2}{1+x}$$

valid for $x \neq -1$, to $x$ defined in (32), we get, for $1/(2k_n) \leq t \leq 1$,

$$\frac{1}{1+(\gamma/\tilde{\sigma})(Q_n(t)-Q_n(1))/\tilde{a}(k_n/n)}$$
$$= t^{\gamma_0}\left[1 - \left(\frac{\gamma}{\tilde{\sigma}} - \gamma_0\right)\frac{1-t^{\gamma_0}}{\gamma_0} - t^{\gamma_0}\frac{\gamma}{\tilde{\sigma}}k_n^{-1/2}Y_n(t)\right.$$
$$\left. + \frac{((\gamma/\tilde{\sigma}-\gamma_0)(1-t^{\gamma_0})/\gamma_0 + t^{\gamma_0}(\gamma/\tilde{\sigma})k_n^{-1/2}Y_n(t))^2}{t^{\gamma_0}(1+(\gamma/\tilde{\sigma})(Q_n(t)-Q_n(1))/\tilde{a}(k_n/n))}\right].$$

Hence the left-hand side of the second equation in (30) equals

$$(33) \quad \begin{aligned} &\int_0^{(2k_n)^{-1}} \frac{1}{1+(\gamma/\tilde{\sigma})(Q_n(t)-Q_n(1))/\tilde{a}(k_n/n)}\,dt \\ &+ \int_{(2k_n)^{-1}}^1 \left(t^{\gamma_0} - \left(\frac{\gamma}{\tilde{\sigma}} - \gamma_0\right)\frac{t^{\gamma_0}-t^{2\gamma_0}}{\gamma_0} - \frac{\gamma}{\tilde{\sigma}}k_n^{-1/2}t^{2\gamma_0}Y_n(t)\right)dt \\ &+ \int_{(2k_n)^{-1}}^1 \frac{((\gamma/\tilde{\sigma}-\gamma_0)(1-t^{\gamma_0})/\gamma_0 + (\gamma/\tilde{\sigma})k_n^{-1/2}t^{\gamma_0}Y_n(t))^2}{1+(\gamma/\tilde{\sigma})(Q_n(t)-Q_n(1))/\tilde{a}(k_n/n)}\,dt. \end{aligned}$$

From (31) it follows easily that the first integral is $o_p(k_n^{-1/2})$. Direct calculations and (19) show that the second integral equals

$$\frac{1}{\gamma_0+1} - \left(\frac{\gamma}{\tilde{\sigma}} - \gamma_0\right)\frac{1}{(\gamma_0+1)(2\gamma_0+1)} - \frac{\gamma}{\tilde{\sigma}}k_n^{-1/2}\int_0^1 t^{2\gamma_0}Y_n(t)\,dt$$
$$+ O_p((2k_n)^{-\gamma_0-1} + k_n^{-1/2}(2k_n)^{-\gamma_0-1} + k_n^{-1/2}(2k_n)^{-\gamma_0-1/2+\varepsilon}).$$

Here, for $\varepsilon < 1/2$, the $O_p$-term is $o_p(k_n^{-1/2})$. By Lemma 3.2, the last integral of (33) is bounded by

$$O_p\left(\int_{(2k_n)^{-1}}^1 t^{\gamma_0}\left(\left(\frac{\gamma}{\tilde{\sigma}} - \gamma_0\right)\frac{1-t^{\gamma_0}}{\gamma_0} + \frac{\gamma}{\tilde{\sigma}}k_n^{-1/2}t^{\gamma_0}Y_n(t)\right)^2 dt\right)$$
$$= O_p(k_n^{-1} + k_n^{-1}(1+(2k_n)^{-\gamma_0+2\varepsilon}) + k_n^{-1}(2k_n)^{-\gamma_0-1/2+\varepsilon}) = o_p(k_n^{-1/2}),$$

if $\varepsilon < 1/4 + \gamma_0/2$. Therefore we have proved

$$\int_0^1 \frac{1}{1+(\gamma/\tilde{\sigma})(Q_n(t)-Q_n(1))/\tilde{a}(k_n/n)}\,dt$$

$$= \frac{1}{\gamma_0 + 1} - \left(\frac{\gamma}{\tilde{\sigma}} - \gamma_0\right)\frac{1}{(\gamma_0+1)(2\gamma_0+1)}$$
$$- \frac{\gamma}{\tilde{\sigma}} k_n^{-1/2} \int_0^1 t^{2\gamma_0} Y_n(t)\,dt + o_p(k_n^{-1/2}).$$

Hence, under the given conditions, system (30) is equivalent to

(34)
$$\gamma_0 + \left(\frac{\gamma}{\tilde{\sigma}} - \gamma_0\right)\frac{1}{\gamma_0+1}$$
$$+ \frac{\gamma}{\tilde{\sigma}} k_n^{-1/2} \int_0^1 t^{\gamma_0} Y_n(t)\,dt + o_p(k_n^{-1/2}) = \gamma,$$
$$\frac{1}{\gamma_0+1} - \left(\frac{\gamma}{\tilde{\sigma}} - \gamma_0\right)\frac{1}{(\gamma_0+1)(2\gamma_0+1)}$$
$$- \frac{\gamma}{\tilde{\sigma}} k_n^{-1/2} \int_0^1 t^{2\gamma_0} Y_n(t)\,dt + o_p(k_n^{-1/2}) = \frac{1}{\gamma+1}.$$

Next we prove that (34) implies (28). First note that, in view of (19) and (21), (34) implies

(35)
$$\gamma_0 + \left(\frac{\gamma}{\tilde{\sigma}} - \gamma_0\right)\frac{1}{\gamma_0+1}$$
$$+ \gamma_0 k_n^{-1/2} \int_0^1 t^{\gamma_0} Y_n(t)\,dt + o_p(k_n^{-1/2}) = \gamma,$$
$$\frac{1}{\gamma_0+1} - \left(\frac{\gamma}{\tilde{\sigma}} - \gamma_0\right)\frac{1}{(\gamma_0+1)(2\gamma_0+1)}$$
$$- \gamma_0 k_n^{-1/2} \int_0^1 t^{2\gamma_0} Y_n(t)\,dt + o_p(k_n^{-1/2}) = \frac{1}{\gamma+1}.$$

The first equation and (21) show that $|\gamma - \gamma_0| = O_p(k_n^{-1/2})$, hence $|\gamma - \gamma_0|^2 = o_p(k_n^{-1/2})$. Therefore $1/(\gamma_0+1) - 1/(\gamma+1) = (\gamma-\gamma_0)/(\gamma_0+1)^2 + o(k_n^{-1/2})$ and so (35) implies

$$\gamma - \gamma_0 - \left(\frac{\gamma}{\tilde{\sigma}} - \gamma_0\right)\frac{1}{\gamma_0+1}$$
$$- k_n^{-1/2}\gamma_0 \int_0^1 t^{\gamma_0} Y_n(t)\,dt + o_p(k_n^{-1/2}) = 0,$$
$$\frac{\gamma-\gamma_0}{(\gamma_0+1)^2} - \left(\frac{\gamma}{\tilde{\sigma}} - \gamma_0\right)\frac{1}{(\gamma_0+1)(2\gamma_0+1)}$$
$$- k_n^{-1/2}\gamma_0 \int_0^1 t^{2\gamma_0} Y_n(t)\,dt + o_p(k_n^{-1/2}) = 0.$$

Now straightforward calculations show that a solution of this linear system in $\gamma - \gamma_0$ and $\gamma/\tilde{\sigma} - \gamma_0$ satisfies (28).

Since conversely a solution of type (28) obviously satisfies the condition (21), it is easily seen that it also solves (34) and thus (20).

CASE $-1/2 < \gamma_0 < 0$. Again, in this case system (20) simplifies to (30). Rewrite the left-hand side of the first equation as

$$\int_0^{s_n} \log\left(1 + \frac{\gamma}{\tilde{\sigma}}\frac{Q_n(t) - Q_n(1)}{\tilde{a}(k_n/n)}\right) dt + \int_{s_n}^1 \log(t^{-\gamma_0}) \, dt$$

$$+ \int_{s_n}^1 \log\left(t^{\gamma_0}\left(1 + \frac{\gamma}{\tilde{\sigma}}\frac{Q_n(t) - Q_n(1)}{\tilde{a}(k_n/n)}\right)\right) dt$$

$$= J_1 + (\gamma_0 + O(s_n|\log s_n|)) + J_2$$

and choose $s_n = k_n^{-\delta}$, with $\delta \in (1/2, (4\varepsilon)^{-1})$ for some $\varepsilon \in (0, 1/2)$.

Now we prove that $J_1$ is negligible. Note that since $t \mapsto Q_n(t)$ is constant when $t \in [0, (2k_n)^{-1}]$, (22) is trivially extended to $t \in [0, 1]$ when $\gamma_0 < 0$. By definition $Q_n(t) - Q_n(1) \geq 0$, for all $t \in [0, 1]$ and $\tilde{a}(k_n/n) > 0$. Since by (21) $P\{\gamma < 0\} \to 1$, Lemma 3.2 implies

$$P\left(\left|\log\left(1 + \frac{\gamma}{\tilde{\sigma}}\frac{Q_n(t) - Q_n(1)}{\tilde{a}(k_n/n)}\right)\right| \leq |\log(C_n t^{-\gamma_0})|, t \in [0, 1]\right) \to 1,$$

as $n \to \infty$, so that $\int_0^{s_n}|\log(C_n t^{-\gamma_0})|\, dt = O_P(s_n|\log s_n|) = o_p(k_n^{-1/2})$ gives $J_1 = o_p(k_n^{-1/2})$.

Next we approximate $J_2$. Check that $0 \leq x - \log(1+x) \leq x^2/[2(1 \wedge (1+x))]$ holds for all $x > -1$. Hence, in view of (17) and Lemma 3.2, choosing $1 + x = t^{\gamma_0}[1 + (\gamma/\tilde{\sigma})(Q_n(t) - Q_n(1))/a(k_n/n)]$, we obtain

$$J_2 = \int_{s_n}^1 \left(\left(\frac{\gamma}{\tilde{\sigma}} - \gamma_0\right)\frac{1-t^{\gamma_0}}{\gamma_0} + \frac{\gamma}{\tilde{\sigma}}k_n^{-1/2}t^{\gamma_0}Y_n(t)\right) dt$$

$$+ O_p\left(\int_{s_n}^1 \left[\left(\frac{\gamma}{\tilde{\sigma}} - \gamma_0\right)\frac{1-t^{\gamma_0}}{\gamma_0} + \frac{\gamma}{\tilde{\sigma}}k_n^{-1/2}t^{\gamma_0}Y_n(t)\right]^2 dt\right)$$

$$= \left(\left(\frac{\gamma}{\tilde{\sigma}} - \gamma_0\right)\frac{1}{\gamma_0+1} + O_p(k_n^{-1/2}s_n^{\gamma_0+1})\right)$$

$$+ \left(\frac{\gamma}{\tilde{\sigma}}k_n^{-1/2}\int_0^1 t^{\gamma_0}Y_n(t)\, dt + O_p\left(k_n^{-1/2}\int_0^{s_n} t^{-1/2-\varepsilon}\, dt\right)\right)$$

$$+ O_p\left(\int_{s_n}^1 \left[\left(\frac{\gamma}{\tilde{\sigma}} - \gamma_0\right)\frac{1-t^{\gamma_0}}{\gamma_0} + \frac{\gamma}{\tilde{\sigma}}k_n^{-1/2}t^{\gamma_0}Y_n(t)\right]^2 dt\right)$$

and from the choice of $s_n$, Corollary 3.1, (21) and $\gamma_0 \in (-1/2, 0)$, it follows that the $O_p$-terms are $o_p(k_n^{-1/2})$. Hence we proved that

$$\int_0^1 \log\left(1 + \frac{\gamma}{\tilde{\sigma}}\frac{Q_n(t) - Q_n(1)}{\tilde{a}(k_n/n)}\right) dt$$

$$= \gamma_0 + \left(\frac{\gamma}{\tilde{\sigma}} - \gamma_0\right)\frac{1}{\gamma_0+1} + \frac{\gamma}{\tilde{\sigma}}k_n^{-1/2}\int_0^1 t^{\gamma_0}Y_n(t)\, dt + o_p(k_n^{-1/2}).$$

Now we turn to the second equation in (30). Use a similar decomposition as in the case $\gamma > 0$ of the left-hand side:

$$\int_0^{s_n} \frac{1}{1+(\gamma/\tilde{\sigma})(Q_n(t)-Q_n(1))/\tilde{a}(k_n/n)} \, dt$$
$$+ \int_{s_n}^1 \left(t^{\gamma_0} - \left(\frac{\gamma}{\tilde{\sigma}} - \gamma_0\right)\frac{t^{\gamma_0}-t^{2\gamma_0}}{\gamma_0} - \frac{\gamma}{\tilde{\sigma}} k_n^{-1/2} t^{2\gamma_0} Y_n(t)\right) dt$$
$$+ \int_{s_n}^1 \frac{((\gamma/\tilde{\sigma}-\gamma_0)(1-t^{\gamma_0})/\gamma_0 + (\gamma/\tilde{\sigma})k_n^{-1/2} t^{\gamma_0} Y_n(t))^2}{1+(\gamma/\tilde{\sigma})(Q_n(t)-Q_n(1))/\tilde{a}(k_n/n)} \, dt$$
$$=: K_1 + K_2 + K_3,$$

with $s_n = k_n^{-\delta}$ for some $\delta \in ((2\gamma_0+2)^{-1}, (4\varepsilon-2\gamma_0)^{-1})$ and $\varepsilon \in (0, \gamma_0+1/2)$. Then by Lemma 3.2, $K_1 = O_p(s_n^{\gamma_0+1}) = o_p(k_n^{-1/2})$. Moreover,

$$K_2 = \frac{1}{\gamma_0+1} - \left(\frac{\gamma}{\tilde{\sigma}} - \gamma_0\right)\frac{1}{(\gamma_0+1)(2\gamma_0+1)} - \frac{\gamma}{\tilde{\sigma}} k_n^{-1/2} \int_0^1 t^{2\gamma_0} Y_n(t) \, dt$$
$$+ O_p(s_n^{\gamma_0+1} + k_n^{-1/2} s_n^{2\gamma_0+1} + k_n^{-1/2} s_n^{\gamma_0+1/2-\varepsilon})$$

and, by the choice of $\delta$, we have that the $O_p$-term is $o_p(k_n^{-1/2})$. Finally from Lemma 3.2 and the definition of $s_n$,

$$K_3 = O_p\left(\int_{s_n}^1 t^{\gamma_0} \left[\left(\frac{\gamma}{\tilde{\sigma}} - \gamma_0\right)\frac{1-t^{\gamma_0}}{\gamma_0} + \frac{\gamma}{\tilde{\sigma}} k_n^{-1/2} t^{\gamma_0} Y_n(t)\right]^2 dt\right)$$
$$= O_p(k_n^{-1}(s_n^{3\gamma_0+1} \vee 1) + k_n^{-1}(s_n^{\gamma_0-2\varepsilon}+1) + k_n^{-1}(s_n^{2\gamma_0+1/2-\varepsilon}+1))$$
$$= o_p(k_n^{-1/2}).$$

Hence, the proof can be concluded by the same arguments as in Case $\gamma_0 > 0$.

CASE $\gamma_0 = 0$. In this case we use (20). Apply (twice) the equality $1/(1+x) = 1 - x + x^2/(1+x)$, the inequality $|x - \log(1+x) - x^2/2 + x^3/3| \leq x^4/[4(1 \wedge (1+x)^4)]$, valid for all $x > -1$, to $x = (\gamma/\tilde{\sigma})(Q_n(t)-Q_n(1))/\tilde{a}(k_n/n)$, and use (23) in Lemma 3.2 to obtain for the left-hand side of the first equa-

tion

$$
\begin{aligned}
\frac{1}{\gamma^2} \int_0^1 & \log\left(1 + \frac{\gamma}{\tilde{\sigma}} \frac{Q_n(t) - Q_n(1)}{\tilde{a}(k_n/n)}\right) \\
& - (1+\gamma) \frac{(\gamma/\tilde{\sigma})(Q_n(t) - Q_n(1))/\tilde{a}(k_n/n)}{1 + (\gamma/\tilde{\sigma})(Q_n(t) - Q_n(1))/\tilde{a}(k_n/n)} \, dt \\
= & \int_0^1 -\frac{1}{\tilde{\sigma}} \frac{Q_n(t) - Q_n(1)}{\tilde{a}(k_n/n)} \, dt \\
& + \int_0^1 \left(\frac{1}{2} + \gamma\right) \frac{1}{\tilde{\sigma}^2} \left(\frac{Q_n(t) - Q_n(1)}{\tilde{a}(k_n/n)}\right)^2 dt \\
& - \int_0^1 \left(\frac{2}{3} + \gamma\right) \frac{\gamma}{\tilde{\sigma}^3} \left(\frac{Q_n(t) - Q_n(1)}{\tilde{a}(k_n/n)}\right)^3 dt \\
& + O_p\left(\int_0^1 \left(\frac{\gamma^2}{\tilde{\sigma}^4} \left[\frac{Q_n(t) - Q_n(1)}{\tilde{a}(k_n/n)}\right]^4 + \frac{\gamma^3}{\tilde{\sigma}^5} \left[\frac{Q_n(t) - Q_n(1)}{\tilde{a}(k_n/n)}\right]^5\right) dt\right).
\end{aligned}
\tag{36}
$$

By (17) the first integral in the right-hand side of the last equation equals $-\tilde{\sigma}^{-1} - \tilde{\sigma}^{-1} k_n^{-1/2} \int_0^1 Y_n(t) \, dt + o_p(k_n^{-1/2})$. For the second integral in the right-hand side of (36) consider

$$
\int_0^{s_n} \left(\frac{1}{2} + \gamma\right) \frac{1}{\tilde{\sigma}^2} \left(\frac{Q_n(t) - Q_n(1)}{\tilde{a}(k_n/n)}\right)^2 dt
$$
$$
+ \int_{s_n}^1 \left(\frac{1}{2} + \gamma\right) \frac{1}{\tilde{\sigma}^2} \left(\frac{Q_n(t) - Q_n(1)}{\tilde{a}(k_n/n)}\right)^2 dt,
$$

with $s_n = k_n^{-\delta}$, $\delta \in (1/2, (4\varepsilon)^{-1})$, $\varepsilon \in (0, 1/2)$. Then the first of these last two integrals is $o_p(k_n^{-1/2})$ by (24). In view of (17), (19) and $|\gamma/\tilde{\sigma}| = O_p(k_n^{-1/2})$, the second integral equals

$$
\frac{1 + 2\gamma}{\tilde{\sigma}^2} - \frac{1}{\tilde{\sigma}^2} k_n^{-1/2} \int_0^1 (\log t) Y_n(t) \, dt + o_p(k_n^{-1/2}).
$$

Using a similar reasoning but with $\delta \in (1/2, 3(4\varepsilon)^{-1} \wedge 4(1 + 6\varepsilon)^{-1})$, $\varepsilon \in (0, 1/2)$, the third integral of (36), equals $-4\gamma \tilde{\sigma}^{-3} + o_p(k_n^{-1/2})$. Finally the $O_p$-term of (36) is clearly $o_p(k_n^{-1/2})$ by (24).

Hence we have that (36) equals

$$
-\frac{1}{\tilde{\sigma}} + \frac{1 + 2\gamma}{\tilde{\sigma}^2} - \frac{4\gamma}{\tilde{\sigma}^3}
$$
$$
- \frac{1}{\tilde{\sigma}} k_n^{-1/2} \int_0^1 Y_n(t) \, dt - \frac{1}{\tilde{\sigma}^2} k_n^{-1/2} \int_0^1 (\log t) Y_n(t) \, dt + o_p(k_n^{-1/2}).
$$

To deal with the left-hand side of the second equation of (20), use again the aforementioned equality for $1/(1+x)$ and (23) in Lemma 3.2 to get

$$
\frac{1+\gamma}{\gamma} \left[\int_0^1 \frac{\gamma}{\tilde{\sigma}} \frac{Q_n(t) - Q_n(1)}{\tilde{a}(k_n/n)} \, dt - \int_0^1 \left(\frac{\gamma}{\tilde{\sigma}} \frac{Q_n(t) - Q_n(1)}{\tilde{a}(k_n/n)}\right)^2 dt \right.
$$

$$+ O_p\bigg(\int_0^1 \bigg[\frac{\gamma}{\tilde{\sigma}}\frac{Q_n(t) - Q_n(1)}{\tilde{a}(k_n/n)}\bigg]^3 dt\bigg)\bigg]$$

$$= (1+\gamma)\bigg[\frac{1}{\tilde{\sigma}} + \frac{1}{\tilde{\sigma}}k_n^{-1/2}\int_0^1 Y_n(t)\,dt - \int_0^{s_n} \gamma\bigg(\frac{1}{\tilde{\sigma}}\frac{Q_n(t) - Q_n(1)}{\tilde{a}(k_n/n)}\bigg)^2 dt$$

$$- \int_{s_n}^1 \gamma\bigg(\frac{1}{\tilde{\sigma}}\frac{Q_n(t) - Q_n(1)}{\tilde{a}(k_n/n)}\bigg)^2 dt + O_p(k_n^{-1}(\log k_n)^3)\bigg],$$

where for the $O_p$-term we used (24) and $\gamma = O_p(k_n^{-1/2})$. Next we consider the second and third integral in the last equality, $L_1$ and $L_2$ say. As for $L_1$, it follows from (24) that

$$L_1 = O_p(s_n\gamma(\log k_n)^2) = o_p(k_n^{-1/2}),$$

if $s_n = k_n^{-\delta}$, $\delta \in (0,1)$. As for $L_2$, from (19) with $\varepsilon \in (0, 1/2)$, we get

$$L_2 = -2\frac{\gamma}{\tilde{\sigma}^2} + O_p(k_n^{-3/2}s_n^{-2\varepsilon} + k_n^{-1}) = -2\frac{\gamma}{\tilde{\sigma}^2} + o_p(k_n^{-1/2}).$$

Hence, we proved that

$$\int_0^1 \bigg(\frac{1}{\gamma} + 1\bigg)\frac{(\gamma/\tilde{\sigma})(Q_n(t) - Q_n(1))/\tilde{a}(k_n/n)}{1 + (\gamma/\tilde{\sigma})(Q_n(t) - Q_n(1))/\tilde{a}(k_n/n)} dt$$

$$= \frac{1+\gamma}{\tilde{\sigma}} - \frac{2\gamma}{\tilde{\sigma}^2} + \frac{1}{\tilde{\sigma}}k_n^{-1/2}\int_0^1 Y_n(t)\,dt + o_p(k_n^{-1/2}).$$

Therefore, under the given conditions, a solution of (20) must satisfy

$$(1-\tilde{\sigma}) + 2\gamma - \frac{4\gamma}{\tilde{\sigma}}$$

$$- \tilde{\sigma}k_n^{-1/2}\int_0^1 Y_n(t)\,dt - k_n^{-1/2}\int_0^1 (\log t)Y_n(t)\,dt + o_p(k_n^{-1/2}) = 0,$$

$$(1-\tilde{\sigma}) + \gamma - \frac{2\gamma}{\tilde{\sigma}} + k_n^{-1/2}\int_0^1 Y_n(t)\,dt + o_p(k_n^{-1/2}) = 0.$$

Next note that the first equation implies $\tilde{\sigma} = 1 + O_p(k_n^{-1/2})$, and so $\gamma/\tilde{\sigma} = \gamma + o_p(k_n^{-1/2})$. Simplifying the above equations, we arrive at (29).

The converse assertion is proved as in Case $\gamma_0 > 0$. $\square$

PROOF OF THEOREM 2.1. Recall that in the previous proofs we used the second order approximation (16). Because of $\tilde{a}(k_n/n)/a(k_n/n) - 1 = o(k_n^{-1/2})$, Proposition 3.1 shows that, under the conditions (7) and (9), any solution $(\hat{\gamma}_n^*, \hat{\sigma}_n^*)$ of the likelihood equations such that $k_n^{1/2}|\hat{\gamma}_n^* - \gamma_0|$ and $k_n^{1/2}|\hat{\sigma}_n^*/a(k_n/n) - 1|$ are stochastically bounded, must satisfy (28) and (29). Hence, in view of Lemma 3.1, $\Phi(k_n/n) \sim \tilde{\Phi}(k_n/n)$ and

$$\frac{\hat{\sigma}_n^*}{a(k_n/n)} - \frac{\hat{\sigma}_n^*}{\tilde{a}(k_n/n)} = \frac{\hat{\sigma}_n^*}{\tilde{a}(k_n/n)} o(k_n^{-1/2}) = o(k_n^{-1/2}),$$

also (10)–(13) hold.

Conversely, according to Proposition 3.1, there exists a solution of (20) satisfying (28) [resp. (29)]. This solution corresponds to a solution of the likelihood equations satisfying (10)–(13). $\square$

PROOF OF COROLLARY 2.1. According to Theorem 2.1, the components of the left-hand side of (15) minus deterministic bias terms converge to certain integrals of a Gaussian process, that is, to normal random variables. If $k_n^{1/2}\Phi(k_n/n) \to \lambda$, the bias term of $k_n^{1/2}(\hat{\gamma}_n - \gamma_0)$ tends to $\lambda((\gamma_0+1)^2/\gamma_0) \times \int_0^1 (t^{\gamma_0} - (2\gamma_0+1)t^{2\gamma_0})\Psi(t)\,dt$. Using (8) the result follows by simple calculations. Similarly, the asymptotic bias of the second component can be derived.

To calculate the variance of the limiting normal random variable corresponding to $k_n^{1/2}(\hat{\gamma}_n - \gamma_0)$, let $X(t) = ((\gamma_0+1)^2/\gamma_0)(t^{\gamma_0} - (2\gamma_0+1)t^{2\gamma_0})(W(1) - t^{-(\gamma_0+1)}W(t))$. Then straightforward calculations show that $\operatorname{var}(\int_0^1 X(t)\,dt) = \int_0^1 \int_0^1 E[X(s)X(t)]\,ds\,dt = (\gamma_0+1)^2$.

Likewise, to obtain the asymptotic covariance of $k_n^{1/2}(\hat{\gamma}_n - \gamma_0)$ with $k_n^{1/2}(\hat{\sigma}_n/a(k_n/n) - 1)$, let $Y(t) = ((\gamma_0+1)/\gamma_0)((\gamma_0+1)(2\gamma_0+1)t^{2\gamma_0} - t^{\gamma_0})(W(1) - t^{-(\gamma_0+1)}W(t))$. Then $\operatorname{cov}(\int_0^1 X(s)\,ds, \int_0^1 Y(t)\,dt) = \int_0^1 \int_0^1 E[X(s)Y(t)]\,ds\,dt = -(1+\gamma_0)$. The limiting variance of the scale estimator is obtained similarly. $\square$

PROOF OF THEOREM 2.2. Integration of the various terms of (18) yields, for $\gamma_0 = 0$,

$$k_n^{1/2}\left(\frac{m_n^{(1)}}{\tilde{a}(k_n/n)} - 1\right)$$
$$= k_n^{1/2}\left(\int_0^1 \frac{Q_n(t) - Q_n(1)}{\tilde{a}(k_n/n)}\,dt + \int_0^1 \log t\,dt\right)$$
$$= \int_0^1 (W_n(1) - t^{-1}W_n(t))\,dt + k_n^{1/2}\tilde{\Phi}\left(\frac{k_n}{n}\right)\int_0^1 \Psi(t)\,dt + o_p(1).$$

Similarly, we obtain

$$k_n^{1/2}\left(\frac{m_n^{(2)}}{(\tilde{a}(k_n/n))^2} - 2\right) = -\int_0^1 2\log t(W_n(1) - t^{-1}W_n(t))\,dt$$
$$+ k_n^{1/2}\Phi\left(\frac{k_n}{n}\right)\int_0^1 2\log t\,\Psi(t)\,dt + o_P(1).$$

Now, using Taylor expansions, straightforward calculations give

$$k_n^{1/2}\hat{\gamma}_* = -\int_0^1 (2+\log t)(W_n(1) - t^{-1}W_n(t))\,dt$$
$$+ k_n^{1/2}\Phi\left(\frac{k_n}{n}\right)\int_0^1 (2+\log t)\Psi(t)\,dt + o_P(1),$$

and hence by (12)

$$k_n^{1/2}(\hat{\gamma}_* - \hat{\gamma}_{\mathrm{MLE}}) \xrightarrow{p} 0.$$

The proof of the second statement is similar. $\square$

PROOF OF REMARK 2.4. Under the stated conditions the following analogue of (18) holds:

$$k_n^{1/2}\left(\frac{\log Q_n(t) - \log Q_n(1)}{\tilde{a}(k_n/n)/F^{\leftarrow}(1-k_n/n)} + \log t\right)$$
$$= W_n(1) - t^{-1}W_n(t) - k_n^{1/2}\Phi^*\left(\frac{k_n}{n}\right)\frac{\log^2 t}{2} + o_p(1)t^{-1/2-\varepsilon},$$

for some function $\Phi^*$ such that

(37) $$\Phi^*\left(\frac{k_n}{n}\right) \sim \frac{\tilde{a}(k_n/n)}{F^{\leftarrow}(1-k_n/n)} = O(k_n^{-1/2});$$

see Draisma, de Haan, Peng and Pereira [(1999), Appendix]. Now the results by de Haan and Stadtmüller (1996) imply that $\tilde{a}(k_n/n)$ tends to a positive constant, while $F^{\leftarrow}(1-k_n/n)$ behaves like a multiple of $\log(n/k_n)$. Hence the bias term is asymptotically negligible if $k_n = o(\log^2 n)$, and the assertion can be concluded by the same reasoning as in the proof of Theorem 2.2. $\square$

H. DREES  
FACULTY OF MATHEMATICS, SPST  
UNIVERSITY OF HAMBURG  
BUNDESSTR. 55  
20146 HAMBURG  
GERMANY  
E-MAIL: drees@math.uni-hamburg.de

A. FERREIRA  
EURANDOM  
EINDHOVEN  
THE NETHERLANDS  
AND  
DEPARTMENT OF MATHEMATICS  
1349-017 LISBON  
TECHNICAL UNIVERSITY LISBON  
PORTUGAL  
E-MAIL: anafh@isa.utl.pt

L. DE HAAN  
SCHOOL OF ECONOMICS  
ERASMUS UNIVERSITY ROTTERDAM  
POSTBUS 1738  
3000 DR ROTTERDAM  
THE NETHERLANDS  
E-MAIL: ldehaan@few.eur.nl